\documentclass[french,12pt]{article}
\usepackage{amssymb}

\usepackage{graphicx}
\usepackage[fleqn]{amsmath}

\newtheorem{theorem}{Theorem}

\newtheorem{corollary}{Corollary}

\newtheorem{definition}{Definition}

\newtheorem{lemma}{Lemma}

\newtheorem{proposition}{Proposition}
\newtheorem{remark}{Remark}

\newenvironment{proof}[1][Proof]{\textbf{#1.} }{\ \rule{0.5em}{0.5em}}
\def\alpha{{\Greekmath 010B}}%
\def\beta{{\Greekmath 010C}}%
\def\gamma{{\Greekmath 010D}}%
\def\delta{{\Greekmath 010E}}%
\def\epsilon{{\Greekmath 010F}}%
\def\zeta{{\Greekmath 0110}}%
\def\eta{{\Greekmath 0111}}%
\def\theta{{\Greekmath 0112}}%
\def\iota{{\Greekmath 0113}}%
\def\kappa{{\Greekmath 0114}}%
\def\lambda{{\Greekmath 0115}}%
\def\mu{{\Greekmath 0116}}%
\def\nu{{\Greekmath 0117}}%
\def\xi{{\Greekmath 0118}}%
\def\pi{{\Greekmath 0119}}%
\def\rho{{\Greekmath 011A}}%
\def\sigma{{\Greekmath 011B}}%
\def\tau{{\Greekmath 011C}}%
\def\upsilon{{\Greekmath 011D}}%
\def\phi{{\Greekmath 011E}}%
\def\chi{{\Greekmath 011F}}%
\def\psi{{\Greekmath 0120}}%
\def\omega{{\Greekmath 0121}}%
\def\varepsilon{{\Greekmath 0122}}%
\def\vartheta{{\Greekmath 0123}}%
\def\varpi{{\Greekmath 0124}}%
\def\varrho{{\Greekmath 0125}}%
\def\varsigma{{\Greekmath 0126}}%
\def\varphi{{\Greekmath 0127}}%

\def\nabla{{\Greekmath 0272}}
\def\FindBoldGroup{%
   {\setbox0=\hbox{$\mathbf{x\global\edef\theboldgroup{\the\mathgroup}}$}}%
}

\def\Greekmath#1#2#3#4{%
    \if@compatibility
        \ifnum\mathgroup=\symbold
           \mathchoice{\mbox{\boldmath$\displaystyle\mathchar"#1#2#3#4$}}%
                      {\mbox{\boldmath$\textstyle\mathchar"#1#2#3#4$}}%
                      {\mbox{\boldmath$\scriptstyle\mathchar"#1#2#3#4$}}%
                      {\mbox{\boldmath$\scriptscriptstyle\mathchar"#1#2#3#4$}}%
        \else
           \mathchar"#1#2#3#4%
        \fi 
    \else 
        \FindBoldGroup
        \ifnum\mathgroup=\theboldgroup 
           \mathchoice{\mbox{\boldmath$\displaystyle\mathchar"#1#2#3#4$}}%
                      {\mbox{\boldmath$\textstyle\mathchar"#1#2#3#4$}}%
                      {\mbox{\boldmath$\scriptstyle\mathchar"#1#2#3#4$}}%
                      {\mbox{\boldmath$\scriptscriptstyle\mathchar"#1#2#3#4$}}%
        \else
           \mathchar"#1#2#3#4%
        \fi     	    
	  \fi}

\newif\ifGreekBold  \GreekBoldfalse
\let\SAVEPBF=\pbf
\def\pbf{\GreekBoldtrue\SAVEPBF}%
\def\@@eqncr{\let\@tempa\relax
    \ifcase\@eqcnt \def\@tempa{& & &}\or \def\@tempa{& &}%
      \else \def\@tempa{&}\fi
     \@tempa
     \if@eqnsw
        \iftag@
           \@taggnum
        \else
           \@eqnnum\stepcounter{equation}%
        \fi
     \fi
     \global\tag@false
     \global\@eqnswtrue
     \global\@eqcnt\z@\cr}

\def\TCItag{\@ifnextchar*{\@TCItagstar}{\@TCItag}}
\def\@TCItag#1{%
    \global\tag@true
    \global\def\@taggnum{(#1)}}
\def\@TCItagstar*#1{%
    \global\tag@true
    \global\def\@taggnum{#1}}
\def\tfrac#1#2{{\textstyle {#1 \over #2}}}%
\def\dfrac#1#2{{\displaystyle \frac{#1}{#2}}}%
\def\QATOP#1#2{{#1 \atop #2}}%
\def\QTATOP#1#2{{\textstyle {#1 \atop #2}}}%
\def\QDATOP#1#2{{\displaystyle {#1 \atop #2}}}%
\def\QTATOPD#1#2#3#4{{\textstyle {#3 \atopwithdelims#1#2 #4}}}%
\def\dint{\mathop{\displaystyle \int}}%
\def\dsum{\mathop{\displaystyle \sum }}%

\begin{document}

\title{Modified Bernstein\vspace{-0.4cm}\vspace{-0.4cm} Polynomials and Jacobi
Polynomials \vspace{-1cm}in $q$-Calculus}
\author{Marie-Madeleine Derriennic \\
{\small IRMAR, INSA, 20 Avenue des Buttes de Co\"{e}smes, CS
14315,35043-RENNES, FRANCE}\vspace{-0.4cm}\vspace{-0.2cm}\\
{\small E-mail adress: mderrien@insa-rennes.fr}}
\date{}
\maketitle

\begin{abstract}
We introduce here a generalization of the modified Bernstein polynomials for
Jacobi weights using the $q$-Bernstein basis proposed by G.M. Phillips to
generalize classical Bernstein Polynomials. The function is evaluated at
points which are in geometric progression in $]0,1[$. Numerous properties of
the modified Bernstein Polynomials are extended to their $q$-analogues:
simultaneous approximation, pointwise convergence even for unbounded
functions, shape-preserving property, Voronovskaya theorem,
self-adjointness. Some properties of the eigenvectors, which are $q$%
-extensions of Jacobi polynomials, are given.
\end{abstract}



\section{Introduction}

\ G.M.Phillips has proposed a generalization of Bernstein polynomials based
on the $q$-integers (cf. \cite{ph0}). We introduce here a q-analogue of the
operators which are often called Bernstein Durrmeyer polynomials and denoted 
$M_{n,1}^{\alpha ,\beta }$ (cf. \cite{der0},\cite{ber}).

In all the paper, we shall assume that $q\in ]0,1[$ and, $\alpha ,\beta >-1$
(part 5 excepted)$.$

For any integer $n$, and a function\ $f$ defined on $]0,1[$ we set 
\begin{equation}
{\LARGE M}_{n,q}^{\alpha ,\beta }f(x)=\sum_{k=0}^{n}f_{n,k,q}^{\alpha ,\beta
}\ b_{n,k,q}(x)  \label{formule0}
\end{equation}
where each $f_{n,k,q}^{\alpha ,\beta }$ is a mean of $f$ defined by Jackson
integrals. The polynomials \ $b_{n,k,q}(x)$ are $q$-analogues of the
Bernstein basis polynomials and are defined by\linebreak $b_{n,k,q}(x)=\left[
\QDATOP{n}{k}\right] _{q}x^{k}(1-x)_{q}^{n-k},$ with $\left[ \QDATOP{n}{k}%
\right] _{q}=\dfrac{\left[ n\right] _{q}!}{\left[ k\right] _{q}!\left[ n-k%
\right] _{q}!}$, ($q$-binomial coefficient), $k=0,\ldots ,n$. They verify $%
\sum_{k=0}^{n}b_{n,k,q}(x)=1$ (cf. \cite{ph0}).

We follow the definitions and notations of (\cite{kac}).

For any real $a,$ $\left[ a\right] _{q}=(1-q^{a)}/(1-q)$, $\Gamma
_{q}(a+1)=(1-q)_{q}^{a}/(1-q)^{a},$\newline
(if $a\in \mathbb{N}$ the $q$-integer $\left[ a\right] _{q}$ is $\left[ a%
\right] _{q}=1+q+\cdots +q^{a-1}$ and $\Gamma _{q}(a+1)=(\left[ a\right]
_{q})!$)$;$\newline
$(1-x)_{q}^{a}=\prod_{j=0}^{\infty }(1-q^{j}x)\left/ \prod_{j=0}^{\infty
}(1-q^{j+a}x)\right. \ $and consequently $(1-x)_{q}^{a+b}=$\linebreak $%
(1-x)_{q}^{a}(1-q^{a}x)_{q}^{b}$ holds for any $b$, $(1-x)_{q}^{m}=%
\prod_{j=0}^{m-1}(1-q^{j}x)$ if $m$ is integer.

The notations will be simplified as much as possible, the superscript $%
\alpha ,\beta $ and the index $q$ when $q$ is fixed,\ will be suppressed in
some proofs.

We introduce the positive bilinear form:\vspace{-0.3cm} 
\begin{equation}
\langle f,g\rangle _{q}^{\alpha ,\beta }=q^{(\alpha +1)(\beta
+1)}(1-q)\sum_{i=0}^{\infty }q^{i}q^{i\alpha }(1-q^{i+1})_{q}^{\beta
}f(q^{i+\beta +1})g(q^{i+\beta +1}),
\end{equation}
whenever it is defined. It can be written under the form of two definite $q$%
-integrals 
\begin{eqnarray*}
\langle f,g\rangle _{q}^{\alpha ,\beta } &=&\int_{0}^{^{q^{\beta
+1}}}t^{\alpha }(1-q^{-\beta }t)_{q}^{\beta }f(t)g(t)d_{q}t \\
\text{and }\langle f,g\rangle _{q}^{\alpha ,\beta } &=&q^{(\alpha +1)(\beta
+1)}\int_{0}^{1}t^{\alpha }(1-qt)_{q}^{\beta }f(q^{\beta +1}t)g(q^{\beta
+1}t)d_{q}t.
\end{eqnarray*}
(the definite $q$-integral of a function $f$ is $%
\int_{0}^{a}f(x)d_{q}x=a(1-q)\sum_{i=0}^{\infty }q^{i}f(q^{i}a)$ (cf. \cite
{kac}))

\begin{definition}
We set in formula (\ref{formule0}) :\hfill 
\begin{equation}
f_{n,k,q}^{\alpha ,\beta }=\dfrac{\langle b_{n,k,q},f\rangle _{q}^{\alpha
,\beta }}{\langle b_{n,k,q},1\rangle _{q}^{\alpha ,\beta }}=\dfrac{%
\int_{0}^{1}t^{k+\alpha }\ (1-qt)_{q}^{n-k+\beta }\ f(q^{\beta +1}t)\ d_{q}t%
}{\int_{0}^{1}t^{k+\alpha }\ (1-qt)_{q}^{n-k+\beta }\ d_{q}t},\text{ }%
k=0,...,n,  \label{formule1}
\end{equation}
\begin{equation}
\text{to define\quad }{\LARGE M}_{n,q}^{\alpha ,\beta }f(x)=\sum_{k=0}^{n}%
\dfrac{\langle b_{n,k,q},f\rangle _{q}^{\alpha ,\beta }}{\langle
b_{n,k,q},1\rangle _{q}^{\alpha ,\beta }}b_{n,k,q}(x).  \label{formule}
\end{equation}
\smallskip 
\end{definition}

We see that the polynomial $M_{n,q}^{\alpha ,\beta }f$ is well defined if
there exists $\gamma \geq 0$ such that $x^{\gamma }f(x)$ is bounded on $%
]0,A]\ $for some $A\in ]0,1]$ and $\alpha >\gamma -1.$ Indeed, $x^{\alpha
}f(x)$ is then $q$-integrable for the weight $w_{q}^{\alpha ,\beta
}(x)=x^{\alpha }(1-qx)_{q}^{\beta }$. We will say, in this case, that $f$
satisfies the condition $C(\alpha )$. Also $\langle f,g\rangle _{q}^{\alpha
,\beta }$ is well defined if the product $fg$ satisfies $C(\alpha ),$
particularly if $f^{2}$ and $g^{2}$ do it.

In many cases, the limit of $M_{n,q}^{\alpha ,\beta }f(x)$ when $q$ tends to 
$1$ is :\newline
$M_{n,1}^{\alpha ,\beta }f(x)=\dsum\limits_{k=0}^{n}\left( \dint_{\hspace{%
-0.2cm}0}^{1}t^{k+\alpha }\ (1-t)^{n-k+\beta }\ f(t)\ dt\left/ \dint_{%
\hspace{-0.2cm}0}^{1}t^{k+\alpha }\ (1-t)^{n-k+\beta }\ dt\right. \right)
b_{n,k}(x)$\newline
with $b_{n,k}(x)=\left( \QATOP{n}{k}\right) x^{k}(1-x)^{n-k}.$

Numerous properties of the operator $M_{n,1}^{\alpha ,\beta }$ will be
extended to $M_{n,q}^{\alpha ,\beta }$ in this paper.\newpage

\section{First properties}

For any $n\in \mathbb{N},$ the operator $M_{n,q}^{\alpha ,\beta }$ has the
following properties.

\textbf{- }It is linear, positive and it preserves the constants so it is a
contraction\newline
($\sup\limits_{x\in \left[ 0,1\right] }\left| M_{n,q}^{\alpha ,\beta
}f(x)\right| \leq \sup\limits_{x\in ]0,1[}\left| f(x)\right| ).$

\noindent \textbf{- }It is self-adjoint: $\langle M_{n,q}^{\alpha ,\beta
}f,g\rangle _{q}^{\alpha ,\beta }=\langle f,M_{n,q}^{\alpha ,\beta }g\rangle
_{q}^{\alpha ,\beta }.$

\noindent \textbf{- }It preserves the degrees of the polynomials of degree $%
\leq n.$

The first properties are consequences of the definition. The last one
follows after the following proposition since $D_{q}x^{p}=\left[ p\right]
x^{p-1}$.

\begin{proposition}
\label{deri}If $f$ verifies the condition $C(\alpha ),$ we have : 
\begin{equation}
\hspace{-1cm}D_{q}M_{n,q}^{\alpha ,\beta }f(x)=\dfrac{\left[ n\right] _{q}}{%
\left[ n+\alpha +\beta +2\right] _{q}}q^{\alpha +\beta +2}M_{n-1,q}^{\alpha
+1,\beta +1}\left( \hspace{-0.1cm}D_{q}f\hspace{-0.1cm}\left( \frac{\cdot }{q%
}\right) \hspace{-0.1cm}\right) (qx),x\in \lbrack 0,1],  \label{der}
\end{equation}
where the $q$-derivative of a function $f$ is $D_{q}f(x)=\dfrac{f(qx)-f(x)}{%
(q-1)x}$ if $x\neq 0$.
\end{proposition}

(When $f^{\prime }$ is continuous on $\left[ 0,1\right] $, the limit of
formula (\ref{der}) is, when $q$ tends to $1,$ $\left( M_{n,1}^{\alpha
,\beta }f\right) ^{\prime }(x)=n\left( n+\alpha +\beta +2\right)
^{-1}M_{n-1,1}^{\alpha +1,\beta +1}\left( f^{\prime }\right) (x))$ (cf. \cite
{derr})).

\begin{proof}
We compute $Db_{n,k}(x)=\left[ n\right]
(b_{n-1,k-1}(qx)/q^{k-1}-b_{n-1,k}(qx)/q^{k})$ if\linebreak $1\leq k\leq n-1$%
\ and $Db_{n,0}(x)=-\left[ n\right] b_{n-1,0}(qx),$ $Db_{n,n}(x)=\left[ n%
\right] b_{n-1,n-1}(qx)/q^{n-1}$ to get\linebreak\ $DM_{n}^{\alpha ,\beta
}f(x)=\left[ n\right] \sum\limits_{k=0}^{n-1}b_{n-1,k}(qx)(f_{n,k+1}^{\alpha
,\beta }-f_{n,k}^{\alpha ,\beta })/q^{k}$.

\noindent We denote $\psi _{n,k}^{\alpha ,\beta }(t)=t^{k+\alpha
}(1-qt)_{q}^{n-k+\beta }$, $k=0,\ldots ,n.$

Recall that the $q$-derivative of $g_{1}g_{2}$ is $%
D_{q}(g_{1}g_{2})(x)=D_{q}g_{1}(x)g_{2}(qx)+g_{1}(x)D_{q}g_{2}(x).$ The $q$%
-Beta functions are $B_{q}(u,v)=\int_{0}^{1}t^{u-1}(1-qt)_{q}^{v-1}d_{q}t=%
\Gamma _{q}(u)\Gamma _{q}(v)/\Gamma _{q}(u+v)$. The function $\psi
_{n-1,k}^{\alpha +1,\beta +1}(\frac{t}{q})f(q^{\beta +1}t),$ $t\in ]0,1[,$
extended by $0\ $in $0$ is continuous at $0.$ Hence we may use a $q$%
-integration by parts to write, for $k=0,\ldots ,n-1$ :\newline
$B_{q}(k+\alpha +2,n-k+\beta +1)\left[ n+\alpha +\beta +2\right]
(f_{n,k+1}^{\alpha ,\beta }-f_{n,k}^{\alpha ,\beta })=$

$-q^{k+\alpha }\int_{0}^{1}(D\psi _{n-1,k}^{\alpha +1,\beta +1})(\frac{t}{q}%
)f(q^{\beta +1}t)d_{q}t=\int_{0}^{1}q^{k+\alpha +\beta +2}\psi
_{n-1,k}^{\alpha +1,\beta +1}(t)(Df)(q^{\beta +1}t)d_{q}t$\newline
$-\left[ \psi _{n-1,k}^{\alpha +1,\beta +1}(\frac{t}{q})f(q^{\beta +1}t)%
\right] _{0}^{1}.$

\noindent Hence $(f_{n,k+1}^{\alpha ,\beta }-f_{n,k}^{\alpha ,\beta })=\frac{%
q^{\alpha +\beta +2+k}}{[n+\alpha +\beta +2]}\frac{\int_{0}^{1}t^{k+\alpha
+1}(1-qt)^{n-k+\beta }(Df)(q^{\beta +1}t)d_{q}t}{B_{q}(k+\alpha +2,n-k+\beta
+1)}$ and\medskip \newline
$DM_{n}^{\alpha ,\beta }f(x)=\dfrac{\left[ n\right] q^{\alpha +\beta +2}}{%
\left[ n+\alpha +\beta +2\right] }\dsum\limits_{k=0}^{n-1}\dfrac{\langle
b_{n-1,k},Df(\frac{\cdot }{q})\rangle _{q}^{\alpha +1,\beta +1}}{%
B_{q}(k+\alpha +2,n-k+\beta +1)}b_{n-1,k}(x).$
\end{proof}

\begin{theorem}
\label{total}The following equality holds for any $x\in \left[ 0,1\right] $ :%
\vspace{-0.3cm} 
\begin{equation}
M_{n,q}^{\alpha ,\beta }f(x)=\sum_{j=0}^{\infty }\Phi _{j,n,q}^{\alpha
,\beta }(x)f(q^{j+\beta +1})\text{ }  \label{total1}
\end{equation}
\begin{eqnarray*}
where,\text{ }\Phi _{j,n,q}^{\alpha ,\beta }(x)
&=&u_{j}\sum_{k=0}^{n}v_{k}b_{n,k,q}(q^{j+\beta +1})b_{n,k,q}(x), \\
u_{j} &=&(1-q)q^{j(\alpha +1)}(1-q^{j+1})_{q}^{\beta },\text{ \ }j\in 
\mathbb{N}, \\
v_{k}^{-1} &=&q^{k(\beta +1)}\left[ \QDATOP{n}{k}\right] _{q}B_{q}(k+\alpha
+1,n-k+\beta +1),\text{ \ }k=0,\ldots ,n.
\end{eqnarray*}
Moreover, for any $r\in \mathbb{N},$ the sequence $\Phi _{r,n,q}^{\alpha
,\beta },\Phi _{r-1,n,q}^{\alpha ,\beta },\ldots ,\Phi _{0,n,q}^{\alpha
,\beta }$ is totally positive, that is to say, the collocation matrix $%
\left( \Phi _{r-j,n,q}^{\alpha ,\beta }(x_{i})\right) _{i=1,\ldots
,m,j=0,\ldots ,r}$ is totally positive for any family $(x_{i}),$ $0\leq
x_{1}<\ldots <x_{m}\leq 1.$
\end{theorem}

\begin{proof}
We set $\Phi _{j}=\Phi _{j,q}^{\alpha ,\beta },b_{n,k,q}=b_{k},c=q^{\beta
+1}.$ The formulae (\ref{total1}) come by writing the definite $q$-integrals 
$\langle b_{k},f\rangle _{q}^{\alpha ,\beta }$ as discrete sums in (\ref
{formule}) and the Beta integrals $\langle b_{k},1\rangle _{q}^{\alpha
,\beta }=\left[ \QDATOP{n}{k}\right] _{q}B_{q}(k+\alpha +1,n-k+\beta +1),$ $%
k=0,\ldots ,n$.

For the total positivity of the $\Phi _{j}$, we have to prove that, for any $%
m\in \mathbb{N}$ and any two families $(x_{i})_{i=1,\ldots
,m},(j_{k})_{k=1,\ldots ,m},$ with $0\leq x_{1}\leq \ldots \leq x_{m}\leq
1,\ 0\leq j_{m}\leq \ldots \leq j_{1},$ the determinant $\det (\Phi
_{j_{l}}(x_{i}))_{_{i=1,\ldots ,m,l=1,\ldots ,m}}$ is non negative. From the
multilinearity of the determinants, there is a basic composition formula for
the discrete sums (cf. \cite{karlin}). We have $\det (\Phi
_{j_{l}}(x_{i}))_{_{i=1,\ldots ,m,l=1,\ldots ,m}}=\prod_{l=1}^{m}u_{j_{l}}E$
where

$E=\det (\sum_{k=0}^{n}v_{k}b_{k}(cq^{j_{l}})b_{k}(x_{i}))_{i=1,\ldots
,m,l=1,\ldots ,m}$

$=\sum_{k_{1}=0}^{n}\ldots \sum_{k_{m}=0}^{n}v_{k_{1}}\ldots v_{k_{m}}\det
(b_{k_{i}}(cq^{j_{l}})b_{k_{i}}(x_{i}))_{i=1,\ldots ,m,l=1,\ldots ,m}$

$=$ $\sum_{k_{1}=0}^{n}\ldots \sum_{k_{m}=0}^{n}v_{k_{1}}\ldots
v_{k_{m}}b_{k_{1}}(x_{1})\ldots b_{k_{m}}(x_{m})\det
(b_{k_{i}}(cq^{j_{l}}))_{i=1,\ldots ,m,l=1,\ldots ,m}$

$=\sum_{0\leq k_{1}\leq \ldots \leq k_{m}\leq n}v_{k_{1}}\ldots
v_{k_{m}}\det (b_{k_{l}}(x_{i}))_{i=1,\ldots ,m,l=1,\ldots ,m}\det
(b_{k_{i}}(cq^{j_{l}}))_{i=1,\ldots ,m,l=1,\ldots ,m}.$

\noindent We know that the $q$-Bernstein basis is totally positive (cf. \cite
{ph1}). Hence we have\newline
$\det (b_{k_{l}}(x_{i}))_{i=1,\ldots ,m,l=1,\ldots ,m}\geq 0\ $and also $%
\det (b_{k_{i}}(cq^{j_{l}}))_{i=1,\ldots ,m,l=1,\ldots ,m}\geq 0$ , since $%
cq^{j_{1}}<cq^{j_{2}}<\ldots <cq^{j_{m}}.$ So $E$ is non negative and the
result follows.
\end{proof}

\begin{corollary}
The number of sign changes of the polynomial $M_{n,q}^{\alpha ,\beta }f$ on $%
]0,1[$ is not greater than the number of sign changes of the function $f$.
\end{corollary}

\begin{proof}
For any $r\in \mathbb{N},$ the sequence $\Phi _{r},\Phi _{r-1},\ldots ,\Phi
_{0}$ is totally positive. We deduce that the number of sign changes of the
polynomial $\sum_{j=0}^{r}\Phi _{j}(x)f(q^{j+\beta +1})$ is not greater than
the number of sign changes of the sequence $f(q^{j+\beta +1}),$ $j=0,\ldots
,r,$ hence not greater than the number of sign changes of the function $f$
in $]0,1[$ (cf. \cite{good}). When $r$ tends to infinity this property is
preserved hence is true for $M_{n}f.$
\end{proof}

\begin{corollary}
Let $f$ be a function satisfying the condition $C(\alpha ).$\hfill 

\begin{enumerate}
\item  If $f$ is increasing (respectively decreasing), then the function $%
M_{n,q}^{\alpha ,\beta }f$ is increasing (respectively decreasing).

\item  If $f$ is convex, then the function $M_{n,q}^{\alpha ,\beta }f$ is
convex.
\end{enumerate}
\end{corollary}

\begin{proof}
\textbf{1) }If $f$ is monotone, for any $s\in \mathbb{R}$ the function $f-s$
has at most one sign change. Hence $M_{n}^{\alpha ,\beta }(f-s)=$ $%
M_{n}^{\alpha ,\beta }f-s$ has at most one sign change and $M_{n}^{\alpha
,\beta }f$ is monotone. If $f$ is increasing, $Df(\frac{\cdot }{q})$ is
positive on $]0,q[.$ Since the operators $M_{n}^{\alpha ,\beta }$ are
positive, we obtain for $x\in \lbrack 0,1],$ $M_{n-1}^{\alpha +1,\beta
+1}\left( Df\left( \frac{\cdot }{q}\right) \right) (qx)\geq 0,$ and, using (%
\ref{der}), $DM_{n}^{\alpha ,\beta }f(x)\geq 0.$ So the function $%
M_{n}^{\alpha ,\beta }f$ is increasing.

\textbf{2) }Let suppose the function $f$ is convex. Since $M_{n}^{\alpha
,\beta }$ preserves the degree of the polynomials, for any real numbers $%
\gamma _{1},\gamma _{2}$ there exist $\delta _{1},\delta _{2}$ and a
function $g$ such that $g(x)=f(x)-\delta _{1}x-\delta _{2}$ and $%
M_{n}^{\alpha ,\beta }f(x)-\gamma _{1}x-\gamma _{2}=M_{n}^{\alpha ,\beta
}g(x).$ The number of sign changes of $g$ being at most two, it is the same
for $M_{n}^{\alpha ,\beta }g$. Hence $M_{n}^{\alpha ,\beta }f$ is convex or
concave. Moreover, if a function $\varphi $ is convex (respectively
concave), $D^{2}\varphi (x)=\frac{q^{3}}{(q-1)^{2}x^{2}}(\varphi (q^{2}x)-%
\left[ 2\right] \varphi (qx)+q\varphi (x))$ is $\geq 0$ (respectively $\leq
0).$ Hence $\medskip $

\noindent $M_{n-2}^{\alpha +2,\beta +2}\left( D^{2}\varphi \left( \frac{%
\cdot }{q}\right) \right) (q^{2}x)\geq 0.$ Using (\ref{der}) two times we
obtain $D^{2}M_{n}^{\alpha ,\beta }\varphi (x)\geq 0$ and $M_{n}^{\alpha
,\beta }\varphi $ is not concave.
\end{proof}

\section{Convergence properties}

\begin{theorem}
\label{w}If $f$ is continuous on $[0,1]$ , 
\begin{equation*}
\left\| M_{n,q}^{\alpha ,\beta }f-f\right\| _{\infty }\leq C_{\alpha ,\beta
}\ \omega \left( f,\frac{1}{\sqrt{\left[ n\right] _{q}}}\right) ,
\end{equation*}
where $\left\| f\right\| _{\infty }$ is the uniform norm of $f$ on $\left[
0,1\right] $ and $\omega (f,.)$ is the usual modulus of continuity of $f,$
the constant $C_{\alpha ,\beta }$ being independent of $n,q,f.$
\end{theorem}

\begin{proof}
As $M_{n}^{\alpha ,\beta }$ is positive, O. Shisha and B. Mond theorem can
be applied. It is sufficient to prove that the order of approximation of $f$
by $M_{n}^{\alpha ,\beta }f$ is $O(\frac{1}{\left[ n\right] })$ for the
functions $f_{i}(x)=x^{i}$, $i=0,1,2.$ We compute the polynomials $%
M_{n}^{\alpha ,\beta }f_{i}$ , $i=1$ and $2,$ with the help of (\ref{der})
by $q$-integrations.

\noindent $\left[ n+\alpha +\beta +2\right] (M_{n}^{\alpha ,\beta
}f_{1}(x)-x)=q^{\beta +1}\left[ \alpha +1\right] -x\left[ \alpha +\beta +2%
\right] $ and

\noindent $\left[ n+\alpha +\beta +2\right] \left[ n+\alpha +\beta +3\right]
(M_{n}^{\alpha ,\beta }f_{2}(x)-x^{2})=$

\noindent $\left[ n\right] \left[ 2\right] x(q^{\alpha +2\beta +3}\left[
\alpha \hspace{-0.1cm}+\hspace{-0.1cm}2\right] (1-x)-\left[ \beta \hspace{%
-0.1cm}+\hspace{-0.1cm}1\right] x)+\left[ a\hspace{-0.1cm}+\hspace{-0.1cm}%
\beta \hspace{-0.1cm}+\hspace{-0.1cm}3\right] \left[ \alpha \hspace{-0.1cm}+%
\hspace{-0.1cm}\beta \hspace{-0.1cm}+\hspace{-0.1cm}2\right] x^{2}+q^{2\beta
+2}\left[ \alpha \hspace{-0.1cm}+\hspace{-0.1cm}2\right] \left[ \alpha 
\hspace{-0.1cm}+\hspace{-0.1cm}1\right] .$

The result follows since $0<q<1$ and $0\leq \left[ a\right] \leq \max (a,1)$
if $a\geq 0.$
\end{proof}

\begin{remark}
\hfill 
\end{remark}

\vspace{-0.5cm}In order to have uniform convergence for all continuous
functions on $\left[ 0,1\right] ,$ it is sufficient to have $%
\lim\limits_{n\rightarrow \infty }M_{n,q}^{\alpha ,\beta }f_{i}=f_{i}$ for $%
i=1,2,$ hence $\lim\limits_{n\rightarrow \infty }1/\left[ n\right] _{q}=0.$
This is realized if and only if $q=q_{n}$ and $\lim\limits_{n\rightarrow
\infty }q_{n}=1$. Indeed, for every $n\in \mathbb{N},$ in both cases $%
q^{n}<1/2$ and $q^{n}\geq 1/2,$ we have $1-q<1/\left[ n\right] _{q}\leq
2\max (1-q,\ln 2/n).$ To maximize the order of approximation by the operator 
$M_{n,q_{n}}^{\alpha ,\beta },$ we are interested to have $\left[ n\right]
_{q_{n}}$ of the same order as $n,$ that is to say to have $\rho n<\left[ n%
\right] _{q_{n}}\leq n,$ for some $\rho >0,$ property which holds with the
following property $S$ for $(q_{n}).$

\begin{definition}
The sequence $(q_{n})_{n\in \mathbb{N}},$ has the property $S$ if and only
if there exists $N\in \mathbb{N}$ and $c>0$ such that, for any $n>N,$ $%
1-q_{n}<c/n.$
\end{definition}

\begin{lemma}
\hfill 

The property $S$ holds if and only if the property $P_{1}$ (respectively $%
P_{2})$ holds where :

$P_{1}\ $is ''There exists $N_{1}\in \mathbb{N}$ and $c_{1}>0$ such that,
for any $n>N_{1},$ $\left[ n\right] _{q_{n}}\geq c_{1}n",$

$P_{2}$ is ''There exists $N_{2}\in \mathbb{N}$ and $c_{2}>0$ such that, for
any $n>N_{2},$ $q_{n}^{n}\geq c_{2}".$
\end{lemma}

\begin{proof}
For any $n\in \mathbb{N},$ the function $\xi (x)=(1-x^{n})/(1-x)$ is
increasing on $[0,1[.$ If $S$ holds, we have, for any $n>N_{1}=N,$ $\left[ n%
\right] _{q_{n}}=\xi (q_{n})\geq \xi (1-c/n)\geq n(1-e^{-c})/c$ and $P_{1}$
follows$.$ If $P_{1}$ holds, we have, for any $n>N=N_{1},$ $1/(1-q_{n})\geq %
\left[ n\right] _{q_{n}}\geq c_{1}n$ and $S$ follows$.$ If $P_{2}$ holds, we
have, for any $n>N=N_{2},$ $n(1-q_{n})\leq -n\ln q_{n}<-\ln c_{2}$ and $S$
follows$.$ If $S$ holds, there exists $N_{2}>N$ such that, if $n>N_{2},$ $%
1-q_{n}<1/2$ hence $q_{n}^{n}>e^{-2n(1-q_{n})}>e^{-2c}$ and $P_{2}$ follows.
\end{proof}

\begin{theorem}
\label{limit}If the function $f$ is continuous at the point $x\in ]0,1[,$
then, 
\begin{equation}
\lim\limits_{n\rightarrow \infty }M_{n,q_{n}}^{\alpha ,\beta }f(x)=f(x)
\end{equation}

\noindent in the two following cases :

\begin{enumerate}
\item  If $f$ is bounded on $[0,1]$ and the sequence $(q_{n})$ is such that $%
\lim\limits_{n\rightarrow \infty }q_{n}=1,$

\item  If there exist real numbers $\alpha ^{\prime },\beta ^{\prime }\geq 0$
and a real $\kappa ^{\prime }>0$ such that, for any $x\in ]0,1[,$ $\left|
x^{\alpha ^{\prime }}(1-x)^{\beta ^{\prime }}f(x)\right| \leq \kappa
^{\prime },$ $\alpha ^{\prime }<\alpha +1,$ $\beta ^{\prime }<\beta +1$ and
the sequence $(q_{n})$ owns the property $S.$
\end{enumerate}
\end{theorem}

\begin{theorem}
\label{vor}If the function $f$ admits a second derivative at the point $x\in
]0,1[$ then, in the cases $1$ and $2$ of theorem \ref{limit}, we have the
Voronovskaya-type limit : 
\begin{equation}
\lim\limits_{n\rightarrow \infty }\left[ n\right] _{q_{n}}(M_{n,q_{n}}^{%
\alpha ,\beta }f(x)-f(x))=\dfrac{d}{dx}\left( x^{\alpha +1}(1-x)^{\beta
+1}f^{\prime }(x)\right) \left/ x^{\alpha }(1-x)^{\beta }\right. .
\label{voron}
\end{equation}
\end{theorem}

(The limit operator is the Jacobi differential operator for the weight $%
x^{\alpha }(1-x)^{\beta })$

For the proofs of theorems \ref{limit} and \ref{vor} we need some
preparation.

\begin{proposition}
\label{tn}

We set, for any $n,m\in \mathbb{N}-\{0\}$ and $x\in \lbrack 0,1],$ $q\in
\lbrack 1/2,1[,$%
\begin{equation}
{\LARGE T}_{n,m,q}(x)=\sum_{k=0}^{n}b_{n,k,q}(x)\dfrac{\dint_{\hspace{-0.2cm}%
0}^{1}t^{k+\alpha }\ (1-qt)_{q}^{n-k+\beta }\ (x-t)^{m}d_{q}t}{\dint_{%
\hspace{-0.2cm}0}^{1}t^{k+\alpha }\ (1-qt)_{q}^{n-k+\beta }d_{q}t}.
\end{equation}
For any $m,$ there exists a constant $K_{m}>0,$ independent of $n$ and $q,$
such that : 
\begin{equation*}
\sup\limits_{x\in \lbrack 0,1]}\left| T_{n,m,q}(x)\right| \leq \left\{ 
\begin{array}{l}
K_{m}/\left[ n\right] _{q}^{m/2}\text{ if }m\text{ is even,} \\ 
K_{m}/\left[ n\right] _{q}^{(m+1)/2}\text{ if }m\text{ is odd.}
\end{array}
\right. 
\end{equation*}
\end{proposition}

To prove this proposition we consider the lemmas \ref{recursive} and \ref
{newton}.

\begin{lemma}
\label{recursive}We set, for any $n,m\in \mathbb{N}$ and $x\in \lbrack 0,1],$%
\newline
$T_{n,m,q}^{1}(x)=\sum\limits_{k=0}^{n}b_{n,k,q}(x)\dfrac{%
\int_{0}^{1}t^{k+\alpha }\ (1-qt)_{q}^{n-k+\beta }\ (x-t)_{q}^{m}d_{q}t}{%
\int_{0}^{1}t^{k+\alpha }\ (1-qt)_{q}^{n-k+\beta }d_{q}t}.\medskip $

The following recursion formula holds for any $q\in \lbrack 1/2,1[$ and $%
m\geq 2$ :

$\left[ n+m+\alpha +\beta +2\right] _{q}q^{-\alpha -2m-1}T_{n,m+1,q}^{1}(x)=$

$(-x(1-x)D_{q}T_{n,m,q}^{1}(x)+T_{n,m,q}^{1}(x)(p_{1,m}(x)+x(1-q)\left[
n+\alpha +\beta \right] _{q}\left[ m+1\right] _{q}q^{1-\alpha -m})$\vspace{%
-0.4cm} 
\begin{equation}
\hspace{-0.2cm}%
+T_{n,m-1,q}^{1}(x)p_{2,m}(x)+T_{n,m-2,q}^{1}(x)p_{3,m}(x)(1-q),
\label{recurs}
\end{equation}
where the polynomials $p_{i,m}(x)$, $i=1,2$ and $3$ are uniformly bounded
with regard to $n$ and $q.$
\end{lemma}

\begin{proof}
1) We set $\psi _{k}(t)=t^{k+\alpha }(1-qt)_{q}^{n-k+\beta }$ and $%
l_{k}(x)=b_{n,k}(x)\left/ \int_{0}^{1}\psi _{k}(t)d_{q}t\right. ,\newline
k=0,\ldots ,n$ and $T_{n,m,q}^{1}=T_{m}^{1}.$ \newline
We compute $x(1-x)D_{q}T_{m}^{1}(x)=x(1-x)\left[ m\right]
\sum_{k=0}^{n}l_{k}(x)\int_{0}^{1}\psi _{k}(t)(x-t)_{q}^{m-1}d_{q}t$\newline
$+\sum_{k=0}^{n}l_{k}(x)\int_{0}^{1}\psi _{k}(t)(qx-t)_{q}^{m}(\left[ k%
\right] -\left[ n\right] x)d_{q}t=A+B.$

We have $A=x(1-x)\left[ m\right] T_{m-1}^{1}(x)$ and \newline
$B=q^{-\alpha }\sum_{k=0}^{n}l_{k}(x)\int_{0}^{1}\left( D_{q}\psi
_{k}\right) (t)\,t(1-qt)(qx-t)_{q}^{m}d_{q}t$\newline
$-q^{1-\alpha -m}\left[ n+\alpha +\beta \right] \sum_{k=0}^{n}l_{k}(x)%
\int_{0}^{1}\psi _{k}(t)(qx-t)_{q}^{m+1}d_{q}t$\newline
$+(x(\left[ n+\alpha +\beta \right] q^{2-\alpha -m}-\left[ n\right] )+\left[
-\alpha \right] )\sum_{k=0}^{n}l_{k}(x)\int_{0}^{1}\psi
_{k}(t)(qx-t)_{q}^{m}d_{q}t$\newline
$=B_{1}+B_{2}+B_{3}.$

We $q$-integrate by parts, setting $\sigma (t)=\left( \frac{t}{q}(1-t)(qx-%
\frac{t}{q})_{q}^{m}\right) $. The $q$-integral in $B_{1}$ becomes $%
\int_{0}^{1}D_{q}\psi _{k}(t)\,t(1-qt)(qx-t)_{q}^{m}d_{q}t=\left[ \psi
_{k}(t)\sigma (t)\right] _{0}^{1}-\int_{0}^{1}\psi _{k}(t)(D_{q}\sigma
)(t)d_{q}t$ for each $k=0,\ldots ,n,$

We expand $\sigma (t)=q^{-2m}(x-\tfrac{t}{q})_{q}^{m+2}+q^{-2m}(\left[ 3%
\right] -q^{m+2})x-q^{m})(x-\tfrac{t}{q})_{q}^{m+1}$\newline
$+q^{-2m+1}(x(q^{m-1}+\left[ m\right] (1-q)q^{m})-x^{2}(1+\left[ 2\right]
q(1-q)\left[ m\right] )(x-\tfrac{t}{q})_{q}^{m}$\newline
$-q^{2m+3}x^{2}\left[ m\right] (1-q)(q^{m-2}-x)(x-\tfrac{t}{q})_{q}^{m-1}.$

We obtain $B_{1}=-q^{-\alpha -2m-1}(\left[ m+2\right] T_{m+1}^{1}(x)-\left[
m+1\right] (\left[ 3\right] x-q^{m+2}x-q^{m})T_{m}^{1}(x)$\newline
$-q^{-\alpha -2m}\left[ m\right] x(q^{m-1}+(1-q)\left[ m\right]
q^{m}-x(1+q(1-q)\left[ 2\right] \left[ m\right] )T_{m-1}^{1}$\newline
$+q^{-\alpha +2-2m}\left[ m-1\right] x^{2}\left[ m\right]
(q^{m-2}-x)(1-q)T_{m-2}^{1}(x).$

\noindent Moreover we have $B_{2}=-q^{1-\alpha -m}\left[ n+\alpha +\beta %
\right] (T_{n,m+1}^{1}(x)-(1-q)\left[ m+1\right] xT_{n,m}^{1}(x)),$

\noindent $B_{3}=(x(q^{n}\left[ \beta -m+2\right] -$ $\left[ 2-\alpha -m%
\right] )+\left[ -\alpha \right] )(T_{n,m}^{1}(x)-(1-q)\left[ m\right]
xT_{n,m-1}^{1}(x)).$
\end{proof}

\begin{lemma}
\label{newton}For any $m\in \mathbb{N}$, $q\in \lbrack 1/2,1[,$ $x\in
\lbrack 0,1]$, the expansion of $(x-t)^{m}$ on the Newton basis at the
points $x/q^{i},$ $i=0,\ldots ,m-1$ is: 
\begin{equation}
(x-t)^{m}=\dsum\limits_{k=1}^{m}d_{m,k}(1-q)^{m-k}(x-t)_{q}^{k},
\label{baseq}
\end{equation}
where the coefficients $d_{m,k},$ verify $\left| d_{m,k}\right| \leq d_{m},$ 
$k=1,\ldots ,m,$ and $d_{m}$ does not depend on $x,t$,$q.$
\end{lemma}

\begin{proof}
For $m=1,$ it is obvious. If for some $m\geq 1,$ the relation (\ref{baseq})
holds, we write $x-t=q^{-k}((x-q^{k}t)-(1-q)\left[ k\right] x)$ for $%
k=1,\ldots ,m$ and we obtain $(x-t)^{m+1}$\newline
$=\dsum\limits_{k=1}^{m+1}d_{m+1,k}(1-q)^{m+1-k}(x-t)_{q}^{k}$ with $%
d_{m+1,k}=q^{-k}(qd_{m,k-1}-\left[ k\right] d_{m,k})$ if $k=1,\ldots ,m$ and 
$d_{m+1,m+1}=q^{-m}d_{m,m}.$ Since $\left| d_{m,k}\right| \leq d_{m}$ we
have $\left| d_{m+1,k}\right| \leq d_{m+1}=2^{-m}(m+1)d_{m},$ $k=1,\ldots
,m+1.$
\end{proof}

\textbf{Proof of the proposition \ref{tn}}

At first we prove that, for any $x,$ $\left| T_{m}^{1}(x)\right| \leq H_{m}/%
\left[ n\right] ^{m/2}$ if $m$ is even (respectively $\leq H_{m}/\left[ n%
\right] ^{(m+1)/2}$ if $m$ is odd), where $H_{m}$ does not depend on $n,q,x.$
We have\linebreak $T_{n,0}^{1}(x)=1$ and the formulae for $M_{n}^{\alpha
,\beta }f_{i},i=1,2$ of the proof of theorem \ref{w} give the result for $%
m=1 $ and $2.$ $\ $The product $\left[ n+\alpha +\beta \right]
(1-q)=1-q^{n+\alpha +\beta }$ is positive and bounded by $\max (1,\left|
1-2^{-(1+\alpha +\beta )}\right| ).$ If the result is true for some $p\geq 2$%
, $p$ odd (respectively even) and any $m\leq p,$ the result for $p+1$
follows from the recursion formula (\ref{recurs}) of lemma \nolinebreak \ref
{recursive}.

Then, we write, for any $n,m\in \mathbb{N}$ and $x\in \lbrack 0,1],$ using
lemma \ref{newton}, and $1-q<1/\left[ n\right] ,$ $\left|
T_{n,m,q}(x)\right| \leq d_{m}\sum_{k=1}^{m}(1-q)^{m-k}\left|
T_{k}^{1}(x)\right| \leq d_{m}(T_{m}^{1}(x)+\sum_{k=1}^{m-1}\left[ n\right]
^{-m+k}H_{k}\left[ n\right] ^{-k/2}$

$\leq d_{m}T_{m}^{1}(x)+\sum_{k=1}^{m-1}H_{k}\left[ n\right] ^{-(m+1)/2})$
and the result follows.$\blacksquare $

Now the following lemma is the key.

\begin{lemma}
Let $(q_{n})$ be a sequence owning the property $S$, $x\in ]0,1[$ and $%
\delta \in ]0,1[,$ $\delta <\min (x,1-x)$. Let $\alpha ,\beta ,\alpha
^{\prime },\beta ^{\prime }$ be real numbers such that $\alpha ^{\prime
},\beta ^{\prime }\geq 0,$ $\alpha >\alpha ^{\prime }-1,$ $\beta >\beta
^{\prime }-1.$ We set $\varphi (t)=t^{-\alpha ^{\prime }}(1-t)^{-\beta
^{\prime }},$ $t\in ]0,1[$ and $I_{x,\delta }(t)=1$ if $\left| t-x\right|
>\delta $, $I_{x,\delta }(t)=0$ elsewhere$.$\newline
The sequence $E_{n}(x,\delta )=\sum\limits_{k=0}^{n}b_{n,k,q_{n}}(x)\dfrac{%
\dint\nolimits_{\hspace{-0.2cm}0}^{1}t^{k+\alpha
}(1-q_{n}t)_{q_{n}}^{n-k+\beta }\varphi (q_{n}^{\beta +1}t)I_{x,\delta
}(t)d_{q_{n}}t}{\dint\nolimits_{\hspace{-0.2cm}0}^{1}t^{k+\alpha
}(1-q_{n}t)_{q_{n}}^{n-k+\beta }d_{q_{n}}t}.$

\noindent\ verifies $\lim\limits_{n\rightarrow \infty }nE_{n}(x,\delta )=0$
for any $x$ and $\delta $ such that $0<\delta <x<1-\delta .$
\end{lemma}

\begin{proof}
Let $\overline{\alpha }$ (respectively $\overline{\beta })$ be the smallest
integer such that $\overline{\alpha }\geq \alpha $ \newline
(respectively $\overline{\beta }\geq \beta )$ and $\tau $ (respectively $%
\tau ^{\prime })$ be a real number such that\newline
$\tau >\dfrac{\overline{\alpha }+\overline{\beta }+2}{\alpha -\alpha
^{\prime }+1}$ (respectively $\tau ^{\prime }>\dfrac{\overline{\alpha }+%
\overline{\beta }+2}{\beta -\beta ^{\prime }+1}).$

For any $k=0,\ldots ,n,$ we have $\int_{0}^{1}t^{k+\alpha
}(1-qt)_{q}^{n-k+\beta }d_{q}t\geq \int_{0}^{1}t^{k+\overline{\alpha }%
}(1-qt)_{q}^{n-k+\overline{\beta }}d_{q}t\linebreak =\frac{\left[ k+%
\overline{\alpha }\right] !\left[ n-k+\overline{\beta }\right] !}{\left[ n+%
\overline{\alpha }+\overline{\beta }+1\right] !}\geq \QTATOPD[ ] {n}{k}%
^{-1}(n+\overline{\alpha }+\overline{\beta }+1)^{-(\overline{\alpha }+%
\overline{\beta }+1)}).$ We set $I_{x,\delta }^{-}(t)=1$ if $0<t<x-\delta $
and $I_{x,\delta }^{+}(t)=1$ if $x+\delta <t<1$, $I_{x,\delta
}^{-}(t)=I_{x,\delta }^{+}(t)=0$ elsewhere.

We split the interval $(0,1)$ introducing $e_{n}=1/n^{\tau },e_{n}^{\prime
}=1-1/n^{\tau ^{\prime }},n\in \mathbb{N},$ and we define, using $\ $again $%
\psi _{n,k,q}^{\alpha ,\beta }(t)=t^{k+\alpha }(1-qt)_{q}^{n-k+\beta }$ and $%
l_{n,k,q}(x)=b_{n,k,q}(x)\left/ \int_{0}^{1}\psi _{n,k,q}^{\alpha ,\beta
}(t)d_{q}t\right. $ of lemma \ref{recursive}, $A_{n}^{1}=%
\sum_{k=0}^{n}l_{n,k,q_{n}}(x)\int\nolimits_{0}^{e_{n}}t^{k+\alpha -\alpha
^{\prime }}(1-q_{n}t)_{q_{n}}^{n-k+\beta }d_{q_{n}}t,$\newline
$A_{n}^{2}=\sum_{k=0}^{n}l_{n,k,q_{n}}(x)\int\nolimits_{e_{n}}^{1}t^{k+%
\alpha -\alpha ^{\prime }}(1-q_{n}t)_{q_{n}}^{n-k+\beta }I_{x,\delta
}^{-}(t)d_{q_{n}}t,$\newline
$A_{n}^{3}=\sum_{k=0}^{n}l_{n,k,q_{n}}(x)\int\nolimits_{0}^{e_{n}^{\prime
}}t^{k+\alpha }(1-q_{n}t)_{q_{n}}^{n-k+\beta }\left/ (1-q_{n}^{\beta
+1}t)^{\beta ^{\prime }}\right. I_{x,\delta }^{+}(t)d_{q_{n}}t,$\newline
$A_{n}^{4}=\sum_{k=0}^{n}l_{n,k,q_{n}}(x)\int\nolimits_{e_{n}^{\prime
}}^{1}t^{k+\alpha }(1-q_{n}t)_{q_{n}}^{n-k+\beta }\left/ (1-q_{n}^{\beta
+1}t)^{\beta ^{\prime }}\right. I_{x,\delta }^{+}(t)d_{q}t.$\newline
If $t>x,$ (respectively if $t<x)$ then $t^{-\alpha ^{\prime }}<x^{-\alpha
^{\prime }}$ (respectively $(1-q_{n}^{\beta +1}t)^{\beta ^{\prime }}$%
\linebreak $>(1-q_{n}^{\beta +1}x)^{\beta ^{\prime }}\geq (1-x)^{\beta
^{\prime }}).$ \newline
Hence, we have $E_{n}(x,\delta )\leq (1/2)^{-(\beta +1)\alpha ^{\prime
}}((1-x)^{-\beta ^{\prime }}(A_{n}^{1}+A_{n}^{2})+x^{-\alpha ^{\prime
}}(A_{n}^{3}+A_{n}^{4}))$ if $q_{n}\geq 1/2$, and it is sufficient to prove $%
\lim\limits_{n\rightarrow \infty }nA_{n}^{i}=0$ for $i=1,2,3,4.$

If $q_{n}^{n}\geq c$ and $e_{n}<1/2,$ we have for $k=0,\ldots ,n,$ $%
\int\nolimits_{0}^{e_{n}}t^{k+\alpha -\alpha ^{\prime
}}(1-q_{n}t)_{q_{n}}^{n-k+\beta }d_{q_{n}}t$\linebreak $=q_{n}^{-k(\beta +1)}%
\QTATOPD[ ] {n}{k}_{q_{n}}^{-1}\int\nolimits_{0}^{e_{n}}b_{n,k,q}(q_{n}^{%
\beta +1}t)t^{\alpha -\alpha ^{\prime }}(1-q_{n}t)_{q_{n}}^{\beta
}d_{q_{n}}t\leq \QTATOPD[ ] {n}{k}_{q_{n}}^{-1}\gamma _{1}e_{n}^{\alpha
-\alpha ^{\prime }+1},$ where $\gamma _{1}$ does not depend on $k,n,x,$
since $q_{n}^{(\beta +1)k}\geq c^{\beta +1},$ $0\leq
b_{n,k,q_{n}}(q_{n}^{\beta +1}t)\leq 1,$ and, $(1-q_{n}t)_{q_{n}}^{\beta
}\leq 1$ if $\beta \geq 0$ and $t\in \lbrack 0,1],$ (respectively $%
(1-q_{n}t)_{q_{n}}^{\beta }\leq (1-e_{n})^{-1}\leq 2$ if $\beta <0$ and $%
t\in \lbrack 0,e_{n}]).$ Hence, we have $A_{n}^{1}\leq \gamma _{1}(n+%
\overline{\alpha }+\overline{\beta }+1)^{\overline{\alpha }+\overline{\beta }%
+1}n^{-\tau (\alpha -\alpha ^{\prime }+1)}.$ The choice of $\tau $ and the
property $S$ on $(q_{n})$ give $\lim\limits_{n\rightarrow \infty
}nA_{n}^{1}=0.$

We choose $m\in \mathbb{N}$ such that $m>\tau \alpha ^{\prime }+1$ and we
write\newline
$A_{n}^{2}\leq n^{\tau \alpha ^{\prime }}\delta
^{-2m}\sum\limits_{k=0}^{n}l_{n,k,q_{n}}(x)\int\nolimits_{e_{n}}^{1}t^{k+%
\alpha }(1-q_{n}t)_{q_{n}}^{n-k+\beta }(x-t)^{2m}d_{q_{n}}t$\newline
$\leq n^{\tau \alpha ^{\prime }}\delta ^{-2m}T_{n,2m,q_{n}}(x)\leq
K_{2m}\delta ^{-2m}n^{\tau \alpha ^{\prime }-m}$, hence $\lim\limits_{n%
\rightarrow \infty }nA_{n}^{2}=0$ by the choice of $m.$

Now we have$,$ if $t<e_{n}^{\prime },$ $(1-q_{n}^{\beta +1}t)^{\beta
^{\prime }}>(1-e_{n}^{\prime })^{\beta ^{\prime }}\geq n^{-\tau ^{\prime
}\beta ^{\prime }},$ hence\newline
$A_{n}^{3}\leq n^{\tau ^{\prime }\beta ^{\prime
}}\sum\limits_{k=0}^{n}l_{n,k,q_{n}}(x)\int\nolimits_{0}^{e_{n}^{\prime
}}t^{k+\alpha }(1-q_{n}t)_{q_{n}}^{n-k+\beta }I_{x,\delta }^{+}(t)d_{q_{n}}t.
$ We choose $m^{\prime }\in \mathbb{N}$ such that $m^{\prime }>\tau ^{\prime
}\beta ^{\prime }+1$ to have $A_{n}^{3}\leq n^{\tau ^{\prime }\beta ^{\prime
}}\delta ^{-2m^{\prime
}}\sum\limits_{k=0}^{n}l_{n,k,q_{n}}(x)\int\nolimits_{0}^{1}t^{k+\alpha
}(1-q_{n}t)_{q_{n}}^{n-k+\beta }(x-t)^{2m^{\prime }}d_{q_{n}}t$

\noindent $\leq n^{\tau ^{\prime }\beta ^{\prime }}\delta ^{-2m^{\prime
}}T_{n,2m^{\prime },q_{n}}\leq K_{2m^{\prime }}\delta ^{-2m^{\prime
}}n^{\tau ^{\prime }\beta ^{\prime }-m^{\prime }},$ hence $%
\lim\limits_{n\rightarrow \infty }nA_{n}^{3}=0$ by the choice of $m^{\prime
}.$

To finish, we prove that $(1-q_{n}^{\beta -\beta ^{\prime
}+1}t)_{q_{n}}^{\beta ^{\prime }}\leq (1-q_{n}^{\beta +1}t)^{\beta ^{\prime
}}$ for any $t\in \lbrack 0,1].$

If $0\leq \beta ^{\prime }<1,$ we use the $q$-binomial formula (cf. \cite
{ask}) and the inequalities\newline
$\left[ -\beta ^{\prime }\right] \leq -\beta ^{\prime }\ $and $\left[ -\beta
^{\prime }+k\right] /\left[ k\right] \geq \left( -\beta ^{\prime }+k\right)
/k$ for any integer $k\geq 1$. In the other cases, if $l$ is the integer
such that $l\leq \beta ^{\prime }<l+1,$ we use the rules of product of $q$%
-binomials to write\newline
$(1-q_{n}^{\beta -\beta ^{\prime }+1}t)_{q_{n}}^{\beta ^{\prime
}}=(1-q_{n}^{\beta -\beta ^{\prime }+1}t)_{q_{n}}^{l}(1-q_{n}^{\beta -\beta
^{\prime }+1+l}t)_{q_{n}}^{\beta ^{\prime }-l}$ and the result follows.

Then, with the same rules, we write, for any $k=0,\ldots ,n$ and $t\in
\lbrack 0,1],$\newline
$(1-q_{n}t)_{q_{n}}^{n-k+\beta }=(1-q_{n}t)_{q_{n}}^{\beta -\beta ^{\prime
}}(1-q_{n}^{\beta -\beta ^{\prime }+1}t)_{q_{n}}^{\beta ^{\prime
}}(1-q_{n}^{\beta +1}t)_{q_{n}}^{n-k}$ and\newline
$(1-q_{n}t)_{q_{n}}^{n-k+\beta }\left/ (1-q_{n}^{\beta +1}t)^{\beta ^{\prime
}}\right. \leq (1-q_{n}t)_{q_{n}}^{\beta -\beta ^{\prime }}(1-q_{n}^{\beta
+1}t)_{q_{n}}^{n-k}.$ We deduce if $q_{n}^{n}\geq c$ and $e_{n}^{\prime
}>1/2,$ $A_{n}^{4}\leq \sum\limits_{k=0}^{n}q_{n}^{-(\beta +1)k}\QTATOPD[ ] {%
n}{k}_{q_{n}}^{-1}l_{n,k,q_{n}}(x)\int\nolimits_{e_{n}^{\prime
}}^{1}t^{\alpha }(1-q_{n}t)_{q_{n}}^{\beta -\beta ^{\prime
}}b_{n,k,q_{n}}(q_{n}^{\beta +1}t)d_{q_{n}}t$

\noindent $\leq \gamma _{2}(n+\overline{\alpha }+\overline{\beta }+1)^{%
\overline{\alpha }+\overline{\beta }+1}(1-e_{n}^{\prime })^{\beta -\beta
^{\prime }+1}$ where $\gamma _{2}$ does not depend on $k,n,x$. The choice of 
$e_{n}^{\prime }$ and $\tau ^{\prime }$ gives $\lim\limits_{n\rightarrow
\infty }nA_{n}^{4}=0.$
\end{proof}

\underline{\textbf{Proof of theorem \ref{limit}}}

Suppose $f$ is continuous at $x\in ]0,1[.$ Let $\varepsilon >0$ be an
arbitrary real number$.$ There exists $\delta ^{\prime }>0$ such that $%
\left| f(x)-f(t)\right| <\varepsilon $ for any $t\in \lbrack 0,1]$ such that 
$\left| x-t\right| <\delta ^{\prime }.$ Let $\delta =\delta ^{\prime }/2$
and $N^{\prime }\in \mathbb{N}$ such that $(1-q_{n}^{\beta +1})x<\delta $
for $n>N^{\prime }.$ Then we have, if $\left| x-t\right| <\delta $ and $%
n>N^{\prime },$ the inequalities $-\delta <-q_{n}^{\beta +1}\delta
<x-q_{n}^{\beta +1}t=q_{n}^{\beta +1}(x-t)+(1-q_{n}^{\beta +1})x<2\delta $
and $\left| f(x)-f(q_{n}^{\beta +1}t)\right| <\varepsilon .$

Hence, we have, if $\left| f\right| $ is bounded by $\kappa ,\left|
f(x)-f(q_{n}^{\beta +1}t)\right| <\varepsilon +2\kappa I_{x,\delta }(t)$
and, in the case 2, $\left| f(x)-f(q_{n}^{\beta +1}t)\right| <\varepsilon
+(\left| f(x)\right| +\kappa ^{\prime }(q_{n}^{\beta +1}t)^{-\alpha ^{\prime
}}(1-q_{n}^{\beta +1}t)^{-\beta ^{\prime }}I_{x,\delta }(t).$

We apply the operator $M_{n,q_{n}}^{\alpha ,\beta }$ at the function $%
h_{x}(t)=f(t)-f(x).$

We have $\left| M_{n,q_{n}}^{\alpha ,\beta }f(x)-f(x)\right| =\left|
M_{n,q_{n}}^{\alpha ,\beta }h_{x}(x)\right| \leq \left( M_{n,q_{n}}^{\alpha
,\beta }\left| h_{x}\right| \right) (x)$

\hfill $\leq \left\{ 
\begin{array}{l}
\varepsilon +2\kappa T_{n,2,q_{n}}(x)/\delta ^{2}\text{ in the case 1, } \\ 
\varepsilon +\left| f(x)\right| T_{n,2,q_{n}}(x)/\delta ^{2}+\kappa ^{\prime
}E_{n}(x,\delta ,q_{n})\text{ in the case 2.}
\end{array}
\right. \medskip $

The second term (respectively and the third term in the case 2) of the right
hand side vanishes when $\left[ n\right] _{q_{n}}$ tends to infinity. Since $%
\lim\limits_{n\rightarrow \infty }1/\left[ n\right] _{q_{n}}=0$ in both
cases (remark 1), the result follows.$\blacksquare $

\underline{\textbf{Proof of theorem \ref{vor}}}

We write Taylor formula at the point $x,$\newline
$f(t)=f(x)+(t-x)f^{\prime }(x)+(t-x)^{2}f^{\prime \prime
}(x)/2+(t-x)^{2}\varepsilon (t-x)$ where $\lim\limits_{u\rightarrow
0}\varepsilon (u)=0.$

We apply the operator $M_{n,q_{n}}^{\alpha ,\beta }$ at the function $f$ of
the variable $t$ to obtain\newline
$M_{n,q_{n}}^{\alpha ,\beta }f(x)-f(x)=-f^{\prime }(x)T_{n,1,q_{n}}(x)+%
\dfrac{f^{\prime \prime }(x)}{2}T_{n,2,q_{n}}(x)+R_{n}(x)$ where\newline
$R_{n}(x)=M_{n,q_{n}}^{\alpha ,\beta }\zeta _{x}(x)$ with $\zeta
_{x}(t)=(t-x)^{2}\varepsilon (t-x).$ We use $\lim\limits_{q\rightarrow 1}%
\left[ a\right] _{q}=a$ for any $a\in \mathbb{R}$ and we verify, with the
help of the formulae of the proof of theorem \ref{w}, that $\lim\limits_{%
\left[ n\right] _{q_{n}}\rightarrow \infty }\left[ n\right]
_{q_{n}}T_{n,1,q_{n}}(x)=(\alpha +\beta +2)x-\alpha -1$ and $\lim\limits_{%
\left[ n\right] _{q_{n}}\rightarrow \infty }\left[ n\right]
_{q_{n}}T_{n,2,q_{n}}(x)=2x(1-x).$ So, to obtain the result, we have to
prove that $\lim\limits_{\left[ n\right] _{q_{n}}\rightarrow \infty }\left[ n%
\right] _{q_{n}}R_{n}(x)=0.$ We proceed in the same manner as in the proof
of theorem \ref{limit}. For any arbitrary $\eta >0,$ we can find $\delta >0$
such that, for $n$ great enough, $\varepsilon (x-t)<\eta $ if $\left|
x-q_{n}^{\beta +1}t\right| <\delta .$

We obtain the inequality $\left| \zeta _{x}(t)\right| \leq \eta
(x-t)^{2}+(\rho _{x}+\left| f(t)\right| )I_{x,\delta }(q^{-(\beta +1)}t)$
for any $t\in ]0,1[,$ where $\rho _{x}$ is independent of $t$ and $\delta .$
We deduce\newline
$\left[ n\right] _{q_{n}}\left| R_{n}(x)\right| \leq \left\{ 
\begin{array}{l}
\left[ n\right] _{q_{n}}(\eta T_{n,2,q_{n}}(x)+(\rho _{x}+\kappa
)T_{n,4,q_{n}}(x)/\delta ^{4})\text{ in the case 1,} \\ 
\left[ n\right] _{q_{n}}(\eta T_{n,2,q_{n}}(x)+\rho
_{x}T_{n,4,q_{n}}(x)/\delta ^{4})+\kappa ^{\prime }nE_{n}(x,\delta ))\text{%
in the case 2.}
\end{array}
\right. $

The right hand side tends to 2$\eta x(1-x)$ when $n$ (hence $\left[ n\right]
_{q_{n}})$ tends to infinity and is as small as wanted.$\blacksquare $

\begin{remark}
\hfill 
\end{remark}

1) We see that the best order of approximation in (\ref{voron}) is in $1/%
\left[ n\right] _{q_{n}}.$ If $1-q_{n}=1/n^{\gamma }$ with $0<\gamma <1,$
then $\lim\limits_{n\rightarrow \infty }\left[ n\right] _{q_{n}}/n^{\gamma
}=1,$ hence $\left[ n\right] _{q_{n}}$ can be replaced by $n^{\gamma }$ in (%
\ref{voron})$.$ If $1-q_{n}=1/n\log n$ or $1/n^{\gamma }$ with $\gamma >1$,
then $\lim\limits_{n\rightarrow \infty }\left[ n\right] _{q_{n}}/n=1,$ $%
\left[ n\right] _{q_{n}}$ can be replaced by $n$ and we refound exactly the
Voronovskaya-limit property of $M_{n,1}^{\alpha ,\beta }(x)$ (case 1)$.$

2) In the case 2, the theorems \ref{limit} and \ref{vor} are valid for $%
M_{n,1}^{\alpha ,\beta },$ if $wf$ is Lebesgue integrable on $[0,1],$ and
this result is new. (In the proof the Jackson integrals have to be replaced
by Lebesgue integrals)

\begin{theorem}
If $f^{\prime }$ is continuous on $[0,1]$ and $q>1/2,$ then 
\begin{equation*}
\hspace{-0.8cm}\left\| D_{q}(M_{n,q}^{\alpha ,\beta }f)-f^{\prime }\right\|
_{\infty }\leq C_{\alpha ,\beta }^{\prime }\left( \omega \left( f^{\prime },%
\frac{1}{\sqrt{\left[ n\right] _{q}}}\right) +\omega \left( f^{\prime
},1-q\right) \right) +\dfrac{\left[ \alpha +\beta +2\right] _{q}}{\left[ n%
\right] _{q}}\left\| f^{\prime }\right\| _{\infty },
\end{equation*}
where $C_{\alpha ,\beta }^{\prime }$ is a constant independent of $n,q,f.$%
\newline
Hence, if $\lim\limits_{n\rightarrow \infty }q_{n}=1,$ then $%
\lim\limits_{n\rightarrow \infty }D_{q_{n}}(M_{n,q_{n}}^{\alpha ,\beta
}f)=f^{\prime }(x)$ uniformly on $\left[ 0,1\right] .$
\end{theorem}

\begin{proof}
We write, using (\ref{der}), for any $x\in \lbrack 0,1],$\newline
$DM_{n}^{\alpha ,\beta }f(x)-f^{\prime }(x)=\tfrac{\left[ n\right] }{\left[
n+\alpha +\beta +2\right] }\left( M_{n-1}^{\alpha +1,\beta +1}\left( \hspace{%
-0.1cm}Df\left( \hspace{-0.1cm}\frac{\cdot }{q}\hspace{-0.1cm}\right) 
\hspace{-0.1cm}\right) (qx)-Df(x)+Df(x)-f^{\prime }(x)\right) $\newline
$+\left( \tfrac{\left[ n\right] }{\left[ n+\alpha +\beta +2\right] }%
-1\right) f^{\prime }(x).$ Since $0<\tfrac{\left[ n\right] }{\left[ n+\alpha
+\beta +2\right] }<1,$ we have $\left| D(M_{n}^{\alpha ,\beta
}f(x))-f^{\prime }(x)\right| $\newline
$\leq \left| M_{n-1}^{\alpha +1,\beta +1}\left( Df\left( \frac{\cdot }{q}%
\right) \right) (qx)-Df(x)\right| +\left| Df(x)-f^{\prime }(x)\right| +%
\tfrac{\left[ \alpha +\beta +2\right] }{\left[ n\right] }\left| f^{\prime
}(x)\right| .$

The theorem (\ref{w}) for the function $Df\left( \frac{\cdot }{q}\right) $
gives\newline
$\left| M_{n-1}^{\alpha +1,\beta +1}\left( Df\left( \frac{\cdot }{q}\right)
\right) (qx)-Df(\frac{\cdot }{q})(qx)\right| \leq C_{\alpha +1,\beta
+1}\omega \left( Df(\frac{\cdot }{q}),\frac{1}{\sqrt{\left[ n-1\right] }}%
\right) .$ Moreover\newline
$\left| Df(x)-f^{\prime }(x)\right| =\left| f^{\prime }(y)-f^{\prime
}(x)\right| $ for some $y$ with $qx<y<x$ hence $\left| y-x\right| <1-q$ and $%
\left| Df(x)-f^{\prime }(x)\right| \leq \omega (f^{\prime },1-q).$ The
modulus of continuity of $Df\left( \frac{\cdot }{q}\right) $ is linked with
the modulus of continuity of $f^{\prime }$. Indeed, for any $y_{i}\in
\lbrack 0,1]\ $and $i=1,2$, there exists $z_{i},$ such that $,$ $%
y_{i}<z_{i}<y_{i}/q$ and $Df(\frac{\cdot }{q})(y_{i})=f^{\prime }(z_{i}).$
As $\left| z_{1}-z_{2}\right| \leq \left| y_{1}-y_{2}\right| /q+(1-q)/q$ we
get $\omega \left( Df(\frac{\cdot }{q}),t\right) =\sup\limits_{\left|
y_{1}-y_{2}\right| <t}\left| Df(\frac{\cdot }{q})(y_{1})-Df(\frac{\cdot }{q}%
)(y_{2})\right| \leq \sup\limits_{\left| y_{1}-y_{2}\right| <t}(\left|
f^{\prime }(z_{1})-f^{\prime }(z_{2})\right| $\newline
$\leq 2(\omega (f^{\prime },t)+\omega (f^{\prime },1-q))$ and the result
follows$.$
\end{proof}

\begin{corollary}
If $f^{\prime }$ is continuous on $[0,1]$ and $1-q_{n}=o(1/n^{4}),$ then 
\begin{equation*}
\lim\limits_{n\rightarrow \infty }\left( M_{n,q_{n}}^{\alpha ,\beta
}f\right) ^{\prime }(x)=f^{\prime }(x)\text{ uniformly on }[0,1].
\end{equation*}
\end{corollary}

\begin{proof}
For any $x\in \lbrack 0,1],$ there exists $u\in (q_{n}x,x)$ such that\newline
$D_{q_{n}}(M_{n,q_{n}}^{\alpha ,\beta }f)(x)=(M_{n,q_{n}}^{\alpha ,\beta
}f)^{\prime }(u)$ and $(x-u)<1-q_{n}.$

Hence $\left| D_{q_{n}}(M_{n,q_{n}}^{\alpha ,\beta
}f)(x)-(M_{n,q_{n}}^{\alpha ,\beta }f)^{\prime }(x)\right| \leq
(1-q_{n})\left| (M_{n,q_{n}}^{\alpha ,\beta }f)^{\prime \prime }(v)\right| $
for some $v$ and $\left\| D_{q_{n}}(M_{n,q_{n}}^{\alpha ,\beta
}f)-(M_{n,q_{n}}^{\alpha ,\beta }f)^{\prime }\right\| _{\infty }\leq
n^{4}(1-q_{n})\left\| M_{n,q_{n}}^{\alpha ,\beta }f\right\| _{\infty }\leq
n^{4}(1-q_{n})\left\| f\right\| _{\infty },$ via Markov inequality.
\end{proof}

\section{self-adjointness properties}

In this part $q\in ]0,1[$ is independent of $n.$

On the space of polynomials $\langle .,.\rangle _{q}^{\alpha ,\beta }$ is an
inner product. Let $\left( P_{_{r,}q}^{\alpha ,\beta }\right) _{r\in \mathbb{%
N}}$ be the sequence of the orthogonal polynomials for $\langle .,.\rangle
_{q}^{\alpha ,\beta }$ such that degree of $P_{_{r,}q}^{\alpha ,\beta }=r$
and $P_{_{r,}q}^{\alpha ,\beta }(0)=\left[ \QDATOP{r+\alpha }{r}\right] _{q}=%
\dfrac{\left[ \alpha +r\right] _{q}\ldots \left[ \alpha +1\right] _{q}}{%
\left[ r\right] _{q}!}$. We set $\nu _{r}=\left( \langle P_{_{r,}q}^{\alpha
,\beta },P_{_{r,}q}^{\alpha ,\beta }\rangle _{q}^{\alpha ,\beta }\right)
^{1/2}.$

We define $U_{q}^{\alpha ,\beta }$ which is a $q$-analogue of the Jacobi
differential operator by: 
\begin{equation}
U_{q}^{\alpha ,\beta }f(x)=D_{q}\left( x^{\alpha +1}(1-q^{-\beta
-1}x)_{q}^{\beta +1}D_{q}f(\frac{x}{q})\right) \left/ x^{\alpha
}(1-q^{-\beta }x)_{q}^{\beta }\right. .  \label{diff}
\end{equation}

\begin{proposition}
The operator $U_{q}^{\alpha ,\beta }$ is self-adjoint for $\langle
.,.\rangle _{q}^{\alpha ,\beta }.$ It preserves the space of polynomials of
degree $r\in \mathbb{N}.$ Consequently, for any $r\in \mathbb{N},$ $%
P_{_{r,}q}^{\alpha ,\beta }$ is eigenvector of $U_{q}^{\alpha ,\beta }$ for
the eigenvalue$\ \mu _{r}^{\alpha ,\beta }=-q^{-\beta -r}\left[ r\right] _{q}%
\left[ r+\alpha +\beta +1\right] _{q}.$
\end{proposition}

\begin{proof}
\hfill

\noindent $U^{\alpha ,\beta }$ is a $q$-differential operator of order $2$.
We compute with $q$-binomial relations, $U^{\alpha ,\beta }f(x)=\left(
-q^{\alpha -\beta }\left[ \beta +1\right] x+\left[ \alpha +1\right]
(1-q^{-\beta -1}x\right) )Df(x)+(1-q^{-\beta -1}x)\frac{x}{q}D^{2}f(\frac{x}{%
q}),$ hence the operator $U$ preserves the degree of polynomials. If $f$ and 
$g$ are polynomials $\langle Uf,g\rangle $ is well defined. We write, since
the $q$-integration by parts is valid, $\langle U^{\alpha ,\beta }f,g\rangle 
$\linebreak $=\left[ x^{\alpha +1}(1-q^{-\beta -1}x)_{q}^{\beta +1}Df(\frac{x%
}{q})g(x)\right] _{0}^{q^{\beta +1}}-\dint_{\hspace{-0.2cm}0}^{q^{\beta
+1}}(qx)^{\alpha +1}(1-q^{-\beta }x)_{q}^{\beta +1}Df(x)Dg(x)d_{q}x.$ and
the first term vanishes. We compute $U^{\alpha ,\beta }f(x)$ for $f(x)=x^{r}$
to obtain\linebreak\ $U^{\alpha ,\beta }f(x)=q^{1-r}\left[ r\right] (\left[
\alpha +r\right] x^{r-1}(1-x)-q^{-\beta -1}\left[ \beta +1\right] x^{r}$
where the coefficient of $x^{r}$ is the eigenvalue $\mu _{r}^{\alpha ,\beta
}.$
\end{proof}

\begin{proposition}
The eigenvectors of the operators $M_{n,q}^{\alpha ,\beta },n\in \mathbb{N}$%
, are the polynomials $P_{_{r,}q}^{\alpha ,\beta },r\in \mathbb{N}\ $and, if 
$f$ satisfies $c(\alpha ),$ we have \newline
$M_{n,q}^{\alpha ,\beta }f=\sum\limits_{r=0}^{n}\lambda _{n,r}^{\alpha
,\beta }\langle f,P_{_{r,}q}^{\alpha ,\beta }\rangle _{q}^{\alpha ,\beta
}P_{_{r,}q}^{\alpha ,\beta }/\nu _{r}^{2}\ $with the eigenvalues

\noindent $\lambda _{n,r}^{\alpha ,\beta }=q^{r(r+\alpha +\beta +1)}\dfrac{%
\left[ n\right] !}{\left[ n-r\right] !}\dfrac{\Gamma _{q}(n+\alpha +\beta +2)%
}{\Gamma _{q}(n+r+\alpha +\beta +2)}$ if $r\leq n,$ $\lambda _{n,r}^{\alpha
,\beta }=0$ otherwise.
\end{proposition}

\begin{proof}
Since $M_{n}$ is self adjoint and preserve the degree of polynomials, the
orthogonal polynomials $P_{r}$ are eigenvectors. The eigenvalue $\lambda
_{n,r}^{\alpha ,\beta }$ is obtained by computing $M_{n}f(x)$ for $%
f(x)=x^{r}.$

We use (\ref{der}) $r$ times to get $D^{r}M_{n}f(x)=\dfrac{q^{r(\alpha
+\beta +r+1)}\left[ r\right] !\left[ n\right] \ldots \left[ n-r+1\right] }{%
\left[ n+\alpha +\beta +2\right] \ldots \left[ n+r+\alpha +\beta +1\right] }%
. $
\end{proof}

\begin{corollary}
\begin{enumerate}
\item  For any $n,m\in \mathbb{N},$ the operators $M_{n,q}^{\alpha ,\beta }$
and $M_{m,q}^{\alpha ,\beta }$ commute on the space of functions satisfying $%
C(\alpha ).$

\item  For any $n\in \mathbb{N},$ the operators $M_{n,q}^{\alpha ,\beta }$
and $U_{q}^{\alpha ,\beta }$ commute on the space of functions $f$ such that 
$f^{\prime }$ is defined in a neighborhood of 0 and is continuous at the
point \nolinebreak 0.
\end{enumerate}
\end{corollary}

\begin{proof}
2) For any $r\in \mathbb{N}$ the $q$-integrals $\langle f,UP_{r}\rangle $
and $\langle Uf,P_{r}\rangle $ are well defined if $f^{\prime }$ is
continuous at the point 0. We go from one to the other by two $q$%
-integrations by parts which are valid because $\lim\limits_{x\rightarrow
0}Df(\frac{x}{q})=f^{\prime }(0).$ Then we write $UM_{n}f=\sum%
\limits_{r=0}^{n}\lambda _{n,r}\langle f,P_{r}\rangle \mu _{r}P_{r}/\nu
_{r}^{2}=\sum\limits_{r=0}^{n}\lambda _{n,r}\langle Uf,P_{r}\rangle
P_{r}/\nu _{r}^{2}=M_{n}Uf.$
\end{proof}

\begin{remark}
\hfill
\end{remark}

This proposition and its corollary open a field to study $\lim_{n\rightarrow
\infty }M_{n,q}^{\alpha ,\beta }f$ for $q$ fixed.\newline
Formally we have $\lim_{n\rightarrow \infty }M_{n,q}^{\alpha ,\beta
}f=S_{q}^{\alpha ,\beta }f=\sum_{r=0}^{\infty }q^{r(r+\alpha +\beta
+1)}\langle f,P_{_{r,}q}^{\alpha ,\beta }\rangle _{q}^{\alpha ,\beta
}P_{_{r,}q}^{\alpha ,\beta }/\nu _{r}^{2}$ and $\lim_{n\rightarrow \infty
}M_{n,q}^{\alpha ,\beta }Q=\sum_{r=0}^{\deg Q}q^{r(r+\alpha +\beta
+1)}\langle Q,P_{_{r,}q}^{\alpha ,\beta }\rangle _{q}^{\alpha ,\beta
}P_{_{r,}q}^{\alpha ,\beta }(x)/\nu _{r}^{2}$ if $Q$ is a polynomial. So $%
\lim_{n\rightarrow \infty }M_{n,q}^{\alpha ,\beta }f=f,$ if and only if $f$
is a constant. Moreover, $\lambda _{n-1,r}^{\alpha ,\beta }-\lambda
_{n,r}^{\alpha ,\beta }$\linebreak $=\dfrac{q^{n+\beta }\lambda
_{n,r}^{\alpha ,\beta }\mu _{r}^{\alpha ,\beta }}{\left[ n\right] \left[
n+\alpha +\beta +1\right] },$ $r\in \mathbb{N},$ hence $M_{n-1}^{\alpha
,\beta }f-M_{n}^{\alpha ,\beta }f=\dfrac{q^{n+\beta }}{\left[ n\right] \left[
n+\alpha +\beta +1\right] }U^{\alpha ,\beta }M_{n}^{\alpha ,\beta }f$ and it
is easy to prove (cf. \cite{ber}) that, when $f^{\prime }$ is defined in a
neighborhood of $0,$ continuous at $0,$ $\left\| M_{n,q}^{\alpha ,\beta
}f-S_{q}^{\alpha ,\beta }f\right\| _{\infty }\leq \gamma _{n}\sup_{x\in %
\left[ 0,q^{\beta +2}\right] }\left| U_{q}^{\alpha ,\beta }f(x)\right| ,$ 
\newline
where $\gamma _{n}=\sum_{k=n+1}^{\infty }\dfrac{q^{k+\beta }}{\left[ k\right]
_{q}\left[ k+\alpha +\beta +1\right] _{q}}\thicksim \dfrac{q^{n+\beta +1}}{%
\left[ n\right] _{q}}$. Of course $U_{q}^{\alpha ,\beta }f$ has to be
bounded on $\left[ 0,q^{\beta +2}\right] ,$ which is true, for example, if $%
f $ is bounded on $\left[ 0,1\right] $ and continuous on $\left[ 0,A\right] $
for some $A<1.$

\begin{proposition}
The polynomials $P_{r,q}^{\alpha ,\beta }$ are $q$-extensions of Jacobi
polynomials for the weight $x^{\alpha }(1-x)^{\beta }$ denoted $%
P_{r}^{\alpha ,\beta },r\in \mathbb{N}.$ They own the following properties
which are the $q$-analogues of the well-known properties of Jacobi
polynomials.

\begin{enumerate}
\item  For any $r\in \mathbb{N},$ $\lim\limits_{q\rightarrow
1}P_{r,q}^{\alpha ,\beta }=$ $P_{r}^{\alpha ,\beta }$,

\item  For any $r\in \mathbb{N},$ the polynomial $P_{r,q}^{\alpha ,\beta }$
is a $q$-hypergeometric function (cf. \cite{kac}) : 
\begin{equation*}
P_{r,q}^{\alpha ,\beta }(x)=\left[ \QTATOP{\alpha +r}{r}\right] {}_{2}\Phi
_{1}\left[ 
\begin{tabular}{ccc}
\begin{tabular}{c}
\begin{tabular}{cc}
$\hspace{-0.4cm}q^{-r},$ & \hspace{-0.3cm}$q^{r+\alpha +\beta +1}$%
\end{tabular}
\\ 
$\hspace{-1.1cm}q^{\alpha +1}$%
\end{tabular}
$\hspace{-0.4cm};$ & \hspace{-0.2cm}$q\,,$ & $q^{-\beta }x\;$%
\end{tabular}
\right]
\end{equation*}
So we have $P_{r,q}^{\alpha ,\beta }(x)=\left[ \QTATOP{\alpha +r}{r}\right]
p_{r}(q^{-\beta -1}x;q^{\alpha +1},q^{\beta +1}:q),$ where $p_{r}(x;u,v:q)$
is the shifted little q-Jacobi polynomial of degree $r$ (cf. \cite{ask},
p.592).

\item  They verify a $q$-analogue of Rodrigues formula : 
\begin{equation*}
P_{r,q}^{\alpha ,\beta }(x)=\dfrac{1}{\left[ r\right] !}\dfrac{%
D_{q}^{r}\left( x^{\alpha +r}(1-q^{-\beta -r}x)_{q}^{\beta +r}\right) }{%
x^{\alpha }(1-q^{-\beta }x)_{q}^{\beta }}.
\end{equation*}

\item  We have the relation for the $q$-derivative : 
\begin{equation*}
D_{q}P_{r,q}^{\alpha ,\beta }\left( \dfrac{\cdot }{q}\right) =-q^{-\beta -r}%
\left[ r+\alpha +\beta +1\right] P_{r-1,q}^{\alpha +1,\beta +1}.
\end{equation*}
\end{enumerate}
\end{proposition}

\begin{proof}
2) We look for the analytic solutions of the equation $U_{q}^{\alpha ,\beta
}f-\mu _{r,q}^{\alpha ,\beta }f=\nolinebreak 0.$\newline
We write $f(x)=\sum_{k=0}^{\infty }a_{k}x^{k}$ and $U_{q}^{\alpha ,\beta
}f(x)-\mu _{r,q}^{\alpha ,\beta }f(x)=\left[ \alpha +1\right] a_{1}-\mu
_{r,q}^{\alpha ,\beta }a_{0}\linebreak +q^{-\beta }\sum_{k=1}^{\infty
}\left( \left[ k+1\right] \left[ k+\alpha +1\right] q^{\beta }a_{k+1}-(\left[
k\right] \left[ k+\alpha +\beta +1\right] -\mu _{r,q}^{\alpha ,\beta
}q^{k+\beta })a_{k}\right) q^{-k}x^{k}$. We obtain $\dfrac{a_{k+1}}{a_{k}}%
=-q^{-r-\beta }\dfrac{(\left[ r\right] -\left[ k\right] )(\left[ k+r+\alpha
+\beta +1\right] }{\left[ k+1\right] \left[ k+\alpha +1\right] }=q^{-\beta
}R(q^{k}),$ for any $k\in \mathbb{N},$ with $R(t)=\dfrac{%
(q^{-r}-t^{-1})(q^{r+\alpha +\beta +1}-t^{-1})}{(q-t^{-1})(q^{\alpha
+1}-t^{-1})}$ and the result follows. \medskip

3) For any polynomial $Q$ of degree $<r$, we verify, with the help of $q$%
-integrations by parts that $\langle \frac{D_{q}^{r}\left( x^{\alpha
+1}(1-q^{-\beta -r}x)_{q}^{\beta +r}\right) }{x^{\alpha }(1-q^{-\beta
}x)_{q}^{\beta }},Q\rangle _{q}^{\alpha ,\beta }=0$.

We compute $P_{r,q}^{\alpha ,\beta }(0)$ by using a $q$-extension of Leibniz
formula. We write\linebreak\ $D_{q}^{r}\left( x^{\alpha +1}(1-q^{-\beta
-r}x)_{q}^{\beta +r}\right) =\sum_{k=0}^{r}\left[ \QATOP{r}{k}\right] \frac{%
\Gamma _{q}(\alpha +r+1)\Gamma _{q}(\beta +r+1)}{\Gamma _{q}(\alpha
+k+1)\Gamma _{q}(\beta +r-k+1)}x^{\alpha +k}(1-q^{-\beta }x)^{\beta
+r-k}q^{c_{k}},$\linebreak where $c_{k}=k(k+\alpha -\beta -r+(k-1)/2)$ and $%
D_{q}^{r}\left( x^{\alpha +1}(1-q^{-\beta -r}x)_{q}^{\beta +r}\right) \left/
x^{\alpha }(1-q^{-\beta }x)_{q}^{\beta }\right. =\sum_{k=0}^{r}\left[ \QATOP{%
r}{k}\right] \frac{\Gamma _{q}(\alpha +r+1)\Gamma _{q}(\beta +r+1)}{\Gamma
_{q}(\alpha +k+1)\Gamma _{q}(\beta +r-k+1)}x^{k}(1-q^{-\beta
-r+k}x)^{r-k}q^{c_{k}}=A(x).$\newline
We obtain $A(0)=\frac{\Gamma _{q}(\alpha +r+1)}{\Gamma _{q}(\alpha +1)}=%
\left[ r\right] !\left[ \QATOP{r+\alpha }{r}\right] $ hence $P_{r,q}^{\alpha
,\beta }(x)=A(x)/\left[ r\right] !.\medskip $

1) We take $\lim\limits_{q\rightarrow 1}A(x)=\sum_{k=0}^{r}\left( \QATOP{r}{k%
}\right) \frac{\Gamma (\alpha +r+1)\Gamma (\beta +r+1)}{\Gamma (\alpha
+k+1)\Gamma (\beta +r-k+1)}x^{k}(1-x)^{r-k}.$ It is $r!P_{r}^{\alpha ,\beta
}(x)$ (Rodrigues formula).\medskip

4) We use (\ref{der}) to prove that $D_{q}P_{r,q}^{\alpha ,\beta }(\frac{%
\cdot }{q})$ is eigenvector of $M_{r-1,q}^{\alpha +1,\beta +1}.$ Hence, it
is equal to $P_{r-1,q}^{\alpha +1,\beta +1}$ up to a constant. We compute $%
D_{q}P_{r,q}^{\alpha ,\beta }(0)=a_{1}=\mu _{r,q}^{\alpha ,\beta }a_{0},$
hence $a_{1}=-\left[ \QATOP{r+\alpha }{r}\right] q^{-\beta -r}\left[ r\right]
\left[ r+\alpha +\beta +1\right] /\left[ \alpha +1\right] =-q^{-\beta -r}%
\left[ r+\alpha +\beta +1\right] \left[ \QATOP{r+\alpha }{r-1}\right] $ and
the result follows.
\end{proof}

\section{The case $\protect\alpha =\protect\beta =-1$}

In this part, we study the operators $M_{n,q}^{\text{-1,-1}}.$ They are
built with $M_{n+1,q}^{0,0}$ as kantorovich operators are built with
Bernstein operators (formula (\ref{kantorovich})).

\begin{definition}
\hfill

The operator $M_{n,q}^{\text{-1,-1}}$ is defined by replacing $\alpha =\beta 
$ by $-1$ in formula (\ref{formule0}). It is : 
\begin{equation*}
M_{n,q}^{\text{-1,-1}}f(x)=\sum_{k=0}^{n}f_{n,k,q}^{\text{-1,-1}}b_{n,k,q}(x)
\end{equation*}
with $f_{n,0,q}^{\text{-1,-1}}=f(0),$ $f_{n,n,q}^{\text{-1,-1}}=f(1)$ and $\ 
$the coefficients $f_{n,k,q}^{\text{-1,-1}},$ for $k=1,\ldots ,n-1,$ are
given by (\ref{formule1}) taking $\alpha =\beta =-1$.\newline
The bilinear form is $\langle f,g\rangle _{q}^{\text{-1,-1}}=\dint_{0}^{^{1}}%
\dfrac{f(t)g(t)}{t(1-t)}d_{q}t.$
\end{definition}

The polynomial $M_{n,q}^{\text{-1,-1}}f$ is well defined for any function $f$
defined on $\left[ 0,1\right] ,$ bounded in a neighborhood of 0 (condition $%
C($-1)). It verifies $M_{n,q}^{\text{-1,-1}}f(0)=f(0)$ and $M_{n,q}^{\text{%
-1,-1}}f(0)=f(1)$, hence it preserves the affine functions.

\begin{proposition}
If the function $f$ is continuous on $\left[ 0,1\right] ,$ then

$\lim\limits_{\alpha \rightarrow \text{-1}}M_{n,q}^{\alpha ,\alpha
}f(x)=M_{n,q}^{\text{-1,-1}}f(x)$ for any $x\in \left[ 0,1\right] .$
\end{proposition}

\begin{proof}
The $q$-binomial coefficients $b_{n,k,q}(x)$ are positive and form a
partition of the unity. Hence it is sufficient to prove that $%
\lim\limits_{\alpha \rightarrow -1}f_{n,k,q}^{\alpha ,\alpha }=f_{n,k,q}^{%
\text{-1,-1}}$ for any $k.$ For $k=1,\ldots ,n-1,$ we compute $%
f_{n,k,q}^{\alpha ,\alpha }-f_{n,k,q}^{\text{-1,-1}}=\frac{%
\int_{0}^{1}t^{k+\alpha }\ (1-qt)_{q}^{n-k+\alpha }\ \left( f(q^{\alpha
+1}t)-f(t)\right) d_{q}t}{B_{q}(k+\alpha +1,n-k+\alpha +1)}$\linebreak $+%
\frac{\int_{0}^{1}t^{k+\alpha }\ (1-qt)_{q}^{n-k+\alpha }f(t)d_{q}t}{%
B_{q}(k+\alpha +1,n-k+\alpha +1)}-\frac{\int_{0}^{1}t^{k-1}\
(1-qt)_{q}^{n-k-1}f(t)d_{q}t}{B_{q}(k,n-k)}=F_{1}+F_{2}+F_{3}.$

We consider $I_{k}=\frac{\int_{0}^{1}t^{k+\alpha }\ (1-qt)_{q}^{n-k+\alpha
}\ \left( f(q^{\alpha +1}t)-\widetilde{f}_{k}(t)\right) d_{q}t}{%
B_{q}(k+\alpha +1,n-k+\alpha +1)}$, with $\widetilde{f}_{0}(t)=f(0),$ $%
\widetilde{f}_{n}(t)$\linebreak $=f(1)$ and $\widetilde{f}_{k}(t)=f(t),$ $%
k=1,\ldots ,n-1$ and prove that $\lim\limits_{\alpha \rightarrow -1}I_{k}=0,$
for $k=0,\ldots ,n.$ We use the additivity of the modulus of continuity of $%
f,$ Beta integrals and we set\linebreak\ $\delta =\left[ \alpha +1\right] /%
\left[ n+2\alpha +2\right] .$ We have $\left| f(q^{\alpha +1}t)-f(0)\right|
\leq \omega (f,t)\leq \omega (f,\delta )(1+t/\delta )$, hence $\left|
I_{0}\right| \leq \omega (f,\delta )(1+\frac{1}{\delta }\int_{0}^{1}t^{%
\alpha +1}\ (1-qt)_{q}^{n+\alpha }d_{q}t\left/ \int_{0}^{1}t^{\alpha }\
(1-qt)_{q}^{n+\alpha }\ d_{q}t\right. )\leq 2\omega (f,\frac{\left[ \alpha +1%
\right] }{\left[ n+2\alpha +2\right] }).$

For $k=n$, we have $\left| f(q^{\alpha +1}t)-f(1)\right| \leq \omega
(f,1-q^{\alpha +1}t),$ hence\newline
$\left| I_{n}\right| \leq \omega (f,\delta )(1+\int_{0}^{1}t^{\alpha +n}\
(1-qt)_{q}^{\alpha +1}d_{q}t\left/ (\delta \int_{0}^{1}t^{\alpha +n}\
(1-qt)_{q}^{\alpha }\ d_{q}t\right. ))\leq 2\omega (f,\frac{\left[ \alpha +1%
\right] }{\left[ n+2\alpha +2\right] }).$ For $k=1,\ldots ,n-1,$ we have $%
\left| f(q^{\alpha +1}t)-f(t)\right| \leq \omega (f,1-q^{\alpha +1})$ and 
\newline
$\left| I_{k}\right| \leq \omega (f,1-q^{\alpha +1}).$ As $f_{n,0,q}^{\alpha
,\alpha }-f(0)=I_{0}$ and $f_{n,n,q}^{\alpha ,\beta }-f(1)=I_{n},$ the
result follows for $k=0$ and $n.$

For the other cases, $F_{1}=I_{k}$ vanishes when $\alpha $ tends to $-1.$
The upper term of $F_{2}$ is the $q$-integral $(1-q)\sum_{j=0}^{n}q^{j(k+%
\alpha +1)}(1-q^{j+1})^{n-k+\alpha }f(q^{j}).$ This serie is uniformly
convergent, hence its limit when $\alpha $ tends to $-1$ is the upper term
of $F_{3}.$ At last the lower term of $F_{2}$ tends to the lower term of $%
A_{3}$ because $\Gamma _{q}$ is continuous.
\end{proof}

Numerous properties shown in the case $\alpha ,\beta >-1$ are still true if $%
\alpha =\beta =-1$. Some of them are given in the following.

\begin{proposition}
If the function $f$ is continuous at the points $0$ and $1,$ and verifies
the condition $C($-1$),$ we have : 
\begin{equation}
D_{q}M_{n,q}^{\text{-1,-1}}f(x)=M_{n-1,q}^{0,0}\left( D_{q}f\left( \frac{%
\cdot }{q}\right) \right) (qx),x\in \lbrack 0,1].  \label{kantorovich}
\end{equation}
\end{proposition}

\begin{proof}
The expressions for $(f_{n,k+1}^{\alpha ,\beta }-f_{n,k}^{\alpha ,\beta })$
in proposition (\ref{deri}) hold if $\alpha =\beta =-1$ and $k=1,\ldots
,n-2. $ For the two other terms we have\newline
$\left[ n-1\right] (f_{n,1}^{\text{-1,-1}}-f(0))=-\left[ t^{n-1}(f(t)-f(1))%
\right] _{0}^{1}+\int_{0}^{1}(1-qt)_{q}^{n-1}D_{q}f(t)d_{q}t$ and\newline
$\left[ n-1\right] (f_{n,n-1}^{\text{-1,-1}}-f(1))=-\left[
(1-t)_{q}^{n-1}(f(t)-f(0))\right] _{0}^{1}+%
\int_{0}^{1}(qt)^{n-1}D_{q}f(t)d_{q}t.$

The first terms vanish since $f$ is continuous at $0$ and $1.$
\end{proof}

\begin{theorem}
\hfill

\begin{enumerate}
\item  $M_{n,q}^{\text{-1,-1}}f(x)=\sum_{j=0}^{\infty }\Phi _{j,n,q}^{\text{%
-1,-1}}(x)f(q^{j}),$ where $\Phi _{j,n,q}^{\text{-1,-1}}$ is defined in
formula (\ref{total1}) with $\alpha =\beta =-1$. The sequence $\Phi
_{r,n,q}^{\text{-1,-1}},\Phi _{r-1,n,q}^{\text{-1,-1}},\ldots ,\Phi
_{0,n,q}^{\text{-1,-1}}$ is totally positive. Consequently the operator $%
M_{n,q}^{\text{-1,-1}}$ diminishes the number of sign changes and preserves
the monotony.

\item  If $f$ is continuous on $[0,1]$, then $\left\| M_{n,q}^{\text{-1,-1}%
}f-f\right\| _{\infty }\leq Cte\ \omega \left( f,\frac{1}{\sqrt{\left[ n%
\right] _{q}}}\right) ,$\newline
(theorem \ref{w}).

\item  If $\lim\limits_{n\rightarrow \infty }q_{n}=1$ and if the function $f$
is bounded on $[0,1],$ then

\begin{enumerate}
\item  $\lim\limits_{n\rightarrow \infty }M_{n,q_{n}}^{\text{-1,-1}%
}f(x)=f(x) $ if $f$ is continuous at the point $x\in ]0,1[,$\newline
(theorem \ref{limit}).

\item  $\lim\limits_{n\rightarrow \infty }\left[ n\right]
_{q_{n}}(M_{n,q_{n}}^{\text{-1,-1}}f(x)-f(x))=f^{\prime \prime }(x)$ if the
function $f$ admits a second derivative at the point $x\in ]0,1[,$ (theorem 
\ref{vor}).
\end{enumerate}
\end{enumerate}
\end{theorem}

\end{document}